\documentclass[12pt]{article}
\usepackage{mathpple}
\usepackage{times}

\usepackage{theorem}  
\usepackage{amsmath}  
\usepackage{amssymb}  
\usepackage{mathrsfs} 
\usepackage{cite}     
\usepackage{psfig}    

\title{\huge\bf $\,$\\[-5ex]
A Lower Bound on the Density of Sphere~Packings 
via Graph Theory\\[2.5ex]}

\author{
   {\bf Michael Krivelevich}\\
   \small Tel Aviv University\vspace*{-1.0ex}\\
   \small Tel Aviv 69978, Israel\vspace*{-1.0ex}\\
   \small\tt krivelev@post.tau.ac.il\\ \\
\and
   {\bf Simon Litsyn}\\
   \small Tel Aviv University\vspace*{-1.0ex}\\
   \small Tel Aviv 69978, Israel\vspace*{-1.0ex}\\
   \small\tt litsyn@eng.tau.ac.il\\ \\
\and
   {\bf Alexander Vardy}\\
   \small University of California, San Diego\vspace*{-1.0ex}\\
   \small 9500 Gilman Drive, La Jolla, CA\,92093, U.S.A.\vspace*{-1.0ex}\\
   \small\tt vardy@kilimanjaro.ucsd.edu\\[3ex]
}

\date{January 28, 2004\vspace{5ex}}

\theoremstyle{plain} 
\theorembodyfont{\normalfont\slshape}

\newtheorem{thm}{Theorem} 
\newenvironment{theorem}
{\begin{thm}\hspace*{-1ex}{\bf.}}{\end{thm}}

\renewcommand{\Bbb}{\mathbb}
\newcommand{\Z}{{\Bbb Z}}

\newcommand{\R}{{\Bbb R}}

\newcommand{\zero}{{\mathbf 0}}

\newcommand{\Tref}[1]{Theorem\,\ref{#1}}

\renewcommand{\le}{\leqslant}

\renewcommand{\ge}{\geqslant}

\newcommand{\deff}{\mbox{$\stackrel{\rm def}{=}$}}

\newcommand{\Strut}[2]{\rule[-#2]{0cm}{#1}}
\newcommand{\highsup}[1]{\raisebox{0.35ex}{\kern 1pt $\scriptstyle {#1} $}}

\newcommand{\shalf}{\mbox{\raisebox{.8mm}{\footnotesize $\scriptstyle 1$}
\footnotesize$\!\!\! / \!\!\!$ \raisebox{-.8mm}{\footnotesize
$\scriptstyle 2$}}}

\newcommand{\sfourth}{\mbox{%
\raisebox{.7mm}{\tiny\hspace{.25ex}$\scriptstyle 1$}%
\tiny$\hspace{-.35ex}/\hspace{-.50ex}$\raisebox{-.7mm}%
{\tiny $\scriptstyle 4$\hspace{.25ex}}}}

\newcommand{\al}{\alpha}

\newcommand{\vol}{{\rm Vol}}

\newcommand{\be}[1]{\begin{equation}\label{#1}}
\newcommand{\ee}{\end{equation}}
\newcommand{\Eq}[1]{(\ref{#1})}

\newcommand{\cG}{{\cal G}}
\newcommand{\cH}{{\cal H}}
\newcommand{\cI}{{\cal I}}

\newcommand{\cP}{{\cal P}}

\newcommand{\cS}{{\cal S}}

\outer\def\proclaim #1. #2\par{\medbreak
 \noindent{\bf#1.\enspace}{\sl#2\par}%
 \ifdim\lastskip<\medskipamount \removelastskip\penalty55\medskip\fi}

\mathchardef\inn="3232
\renewcommand{\in}{\mbox{$\,\inn\,$}}

\makeatletter

\gdef\@punct{.\ \ }  
\def\@sect#1#2#3#4#5#6[#7]#8{%
  \ifnum #2>\c@secnumdepth
     \def\@svsec{}
  \else
     \refstepcounter{#1}\edef\@svsec{%
     \ifnum #2>0{{\csname the#1\endcsname}}.\fi%
    \hskip .5em}
  \fi
  \@tempskipa #5\relax
  \ifdim \@tempskipa>\z@
     \begingroup #6\relax
       \@hangfrom{\hskip #3\relax\@svsec}{\interlinepenalty \@M #8\par}
     \endgroup
     \csname #1mark\endcsname{#7}
     \addcontentsline{toc}{#1}{\ifnum #2>\c@secnumdepth\else
          \protect\numberline{\csname the#1\endcsname}\fi#7}
  \else
     \def\@svsechd{#6\hskip #3\@svsec #8\@punct\csname #1mark\endcsname{#7}
     \addcontentsline{toc}{#1}{\ifnum #2>\c@secnumdepth \else
          \protect\numberline{\csname the#1\endcsname}\fi#7}}
  \fi
  \@xsect{#5}}

\def\@ssect#1#2#3#4#5{\@tempskipa #3\relax
  \ifdim \@tempskipa>\z@
    \begingroup #4\@hangfrom{\hskip #1}{\interlinepenalty \@M #5\par}\endgroup
  \else \def\@svsechd{#4\hskip #1\relax #5\@punct}\fi
  \@xsect{#3}}

\makeatother

\begin{document}

\maketitle
\thispagestyle{empty}

\begin{center}
{\bf Abstract} \vspace{-5.50ex}
\end{center}
\begin{list}{}
{
\addtolength{\leftmargin}{-4.00ex}
\setlength{\rightmargin}{\leftmargin}
}
\item\noindent

\looseness=-1
Using graph-theoretic methods we give a new proof that for all sufficiently
large~$n$,~there exist sphere packings in $\R^n$ of density at least
$cn2^{-n}$ exceeding the classical Minkowski bound by
a factor linear in $n$. This matches up to a constant the
best known lower bounds on the density of sphere packings due to 
Rogers~\cite{Rogers}, Davenport-Rogers~\cite{DR}, and Ball~\cite{Ball}. 
The suggested method makes it possible to describe the points of such
a~packing with complexity $\exp(n\log n)$, 
which is significantly lower than in the other approaches.
\end{list}
\vspace{4ex}

\vfill


\footnotetext[1]
{\\
This work was supported in part by the USA -- Israel Binational Science 
Foundation (M.K.), the Israel Science Foundation (M.K., S.L.),  
the Packard Foundation and the National Science Foundation (A.V.).\vspace*{-3ex}
}

 
\thispagestyle{empty}
\addtocounter{page}{-1}
\newpage


\section{Introduction}
\vspace{-1.5ex} 
\label{sec1}

A sphere packing $\cP$ in $\R^n$ is a collection of non-intersecting
open spheres of equal radii, and its density $\Delta(\cP)$ is the
fraction of space covered by their interiors. Let 
$$
\Delta_n 
\ \ \deff \ \
\sup_{\cP \subset \R^n} \Delta(\cP)
$$
where the supremum is taken over all packings in $\R^n$.
%
\looseness=-1
A celebrated theorem of Minkowski states that 
$\Delta_n \ge \zeta(n)/2^{n-1}$ for all $n \ge 2$.
Since $\zeta(n) = 1 + o(1)$, the asymptotic behavior of
the Minkowski bound~\cite{Minkowski} is given by $\Omega(2^{-n})$.
Asymptotic improvements of the Minkowski bound
were obtained by Rogers~\cite{Rogers}, Davenport and Rogers~\cite{DR}, 
and Ball~\cite{Ball}, all of them being of the form  
$\Delta_n \ge cn 2^{-n}$ 
where $c>0$ is an absolute constant. 
The best currently known lower bound on $\Delta_n$ is due
to Ball~\cite{Ball}, who showed that there exist lattice
packings with density at least $2(n\,{-}\,1)2^{-n}\zeta(n)$.
In this note, we use results from graph theory 
to prove the following theorem.

\begin{theorem}
\label{Thm1}
For all sufficiently large $n$, there exists a sphere 
packing $\cP_n \subset \R^n$ such that
\be{thm1}
\Delta(\cP_n) \ \ge \ 0.01 n 2^{-n}
\ee
Moreover, the spheres in $\cP_n$ can be described 
using a deterministic procedure whose complexity 
is at most $O(2^{\gamma n\log_2\! n})$ for an absolute constant 
$\gamma$.\vspace{-1ex}
\end{theorem}

Although the constant $c = 0.01$ 
in~\Eq{thm1} is not as high 
as in the bounds of Rogers $(c = 0.74)$, Daven\-port-Rogers ($c \,{=}\, 1.68$),
and Ball $(c \,{=}\, 2)$, \Tref{Thm1} still provides an improvement
upon the Minkowski bound by the same linear in $n$ factor. 
With some effort, the constant in \Tref{Thm1} can be increased 
by at least a factor of 10. However, the main merit of our proof 
is not so much in the result itself, but rather in the approach we use. 
Our argument is essentially different from all those previously employed and
is technically very simple. Instead of relying directly on the 
geometry of numbers, we apply graph theoretic tools: specifically, 
we use lower bounds on the independence number of locally sparse
graphs. A similar approach has been recently used
by Jiang and Vardy~\cite{JV} in an asymptotic improvement of 
the classical Gilbert-Varshamov bound in coding theory.

\section{Graph-theoretic proof of the Minkowski bound}
\vspace{-1.5ex} 
\label{sec2}

\looseness=-1
First, let us define two cubes in $\R^n$ --- a smaller cube $K_0$ of side
$s_n$ and a larger cube $K_1$ of side $s_n + 2r_n$. Specifically
\be{cubes}
K_0 \ \ \deff \ \,
\left\{ x \in \R^n ~:~  |x_i| \le \frac{s_n}{2} \Strut{3.00ex}{0ex}\right\}
\hspace{3ex}\text{and}\hspace{3ex} 
K_1 \ \ \deff \ \,
\left\{ x \in \R^n ~:~ |x_i| \le \frac{s_n}{2} + r_n\Strut{3.00ex}{0ex}\right\}
\vspace{.5ex}
\ee
where $r_n$ and $s_n$ are, so far, arbitrary functions of $n$,
except that we assume that $r_n, s_n \in 2\Z$.
Next, define a graph $\cG_n$ as follows: the vertices
of $\cG_n$ are given by $V(\cG_n) = \Z^n \cap K_0$, and
$\{u,v\} \in E(\cG_n)$ iff $d(u,v) < 2r_n$, where $d(\cdot,\cdot)$
is the Euclidean distance~in~$\R^n$.\pagebreak[3.99]
Let $\cI$ be~a~\emph{maximal\/} 
independent set in $\cG_n$ --- that is, $\cI$ is such that every vertex of 
$V(\cG_n)\setminus\cI$ is adjacent to at least one vertex in $\cI$.
By the definition of $E(\cG_n)$, spheres of radius $r_n$ about 
the points of~$\cI$ do not overlap. Moreover, all such spheres 
lie inside the cube $K_1$.
Since $K_1$ tiles $\R^n$, we conclude that there exists
a sphere packing $\cP_n$ 
with density
\be{first-bound}
\Delta (\cP_n)
\ = \
\frac{|\cI|\, V_n (r_n)^n}{\vol(K_1)}
\ = \
\frac{|\cI|\, V_n (r_n)^n}{(s_n + 2r_n)^n}
\ee
where $V_n$ is the volume of a unit sphere in $\R^n$. Let $d_n$
denote the maximum degree of a~vertex in $\cG_n$. It is well-known
(and obvious, for any graph $G$) that
\be{trivial-bound}
|\cI|
\ \ge \
\frac{|V(\cG_n)|}{d_n+1}
\ = \
\frac{(s_n+1)^n}{d_n+1}
\ee
Let $\cS_n(r)$ denote the open sphere of radius $r$ about the origin
$\zero$ of $\R^n$.
If the ratio $s_n/r_n$ is sufficiently large, as we shall assume,
then $d_n+1$ is just the number of points of $\Z^n$ contained in
$\cS_n(2r_n)$.
Thus we can roughly estimate $d_n+1$ as simply the volume
$V_n 2^n (r_n)^n$ of $\cS_n(2r_n)$.
Combining this estimate with~\Eq{first-bound}, \Eq{trivial-bound},
and taking $s_n = 2 n^2 r_n$, we obtain
$$
\Delta(\cP_n)
\,\ \gtrsim \,\
\frac{(s_n)^n V_n (r_n)^n}{ V_n 2^n (r_n)^n (s_n + 2r_n)^n}
\ = \
\frac{1}{2^n \left(1 + \frac{2r_n}{s_n}\right)^{\!n}}
\ \ge \
\frac{1+o(1)}{2^n}
$$
Alternatively, a precise bound on $d_n$ can be derived
as follows. With each point~$v \in \Z^n$, we associate the unit cube
$$
K(v) \ \ \deff \
\left\{ x \in \R^n ~:~ -\frac{1}{2} < x_i - v_i \le \frac{1}{2} \,
\Strut{2.75ex}{0ex}\right\}
$$
Such cubes are fundamental domains of $\Z^n$; hence, they do not intersect.
The length of the main diagonal of $K(v)$ is $\sqrt{n}$, which
implies that $d(x,v) \le \sqrt{n}/2$ for all\pagebreak[3.99]
$x \in K(v)$.
Let $D_n = \cS_n(2r_n) \cap \Z^n$, so that $d_n + 1 = |D_n|$.
It follows by the triangle inequality that if $v \in D_n$, then
$d(x,\zero) < 2r_n + \sqrt{n}/2$ for all $x \in K(v)$.
Hence
\be{d_n-upper}
\bigcup_{v \inn D_n} \hspace{-.75ex} K(v)
\ \subset \
\cS_n \!\left(2r_n + \sqrt{n}/2\right)
\ee
Expressing the volume of $\cup_{v \inn D_n} K(v)$
as $|D_n| = d_n+1$, this implies that the maximum 
degree of a vertex in $\cG_n$ is bounded by
\be{d_n-bounds}
d_n + 1
\ \le \
V_n \left(2r_n + \sqrt{n}/2\right)^n
\ee
Combining \Eq{d_n-bounds} with \Eq{first-bound} and \Eq{trivial-bound}, 
then taking $r_n = 2n^2$ and $s_n = 2n^4$ (say), proves
that the density of $\cP_n$ is at least
\be{Minkowski-bound}
\hspace*{-1.25ex}\Delta(\cP_n)
\: \ge \:
\frac{(s_n)^n V_n (r_n)^n}
{(s_n \hspace{.75pt}{+}\hspace{.75pt} 2r_n)^n V_n 
\left(2r_n \hspace{.75pt}{+}\hspace{.75pt} \sqrt{n}/2\right)^n}
\, = \,
\frac{1}{2^n
\left(1 \,{+}\, \frac{2r_n}{s_n}\right)^{\!n} \!\!
\left(1 \,{+}\, \frac{\sqrt{n}}{4r_n}\right)^{\!n}}
\ = \
\frac{1\,{+}\,o(1)}{2^n}
\ee
Asymptotically, \Eq{Minkowski-bound}
coincides with the Minkowski bound on $\Delta_n$. 
Since a maximal independent set in~$\cG_n$ can be found
in time $O\!\left(|V(\cG_n)|^2\right)$~and
$
|V(\cG_n)|
=
(s_n+1)^n
=
\left(2n^4+1\right)^n
$,
the bound in~\Eq{Minkowski-bound} also reproduces the result 
of Litsyn~and~Tsfa\-sman~\cite[Theorem\,1]{LT}.
\pagebreak[3.99]

\section{Asymptotic improvement using locally sparse graphs}
\vspace{-1.5ex} 
\label{sec3}

\looseness=-1
We next take $\cI$ to be an independent set of \emph{maximum size\/} 
in $\cG_n$ and consider bounds on $|\cI| = \al(\cG_n)$ that are sharper than
the trivial bound in~\Eq{trivial-bound}. Let $T_n$ denote the
number of triangles in $\cG_n$. Then it is
known~\cite[Lemma\,15, p.\,296]{AKS,Bollobas} that
\be{Bollobas}
\alpha(\cG_n)
\ \ge \
\frac{|V(\cG_n)|}{10 d_n}
\left( \Strut{4ex}{0ex}
\log_2 d_n \ - \
\shalf\hspace{0.25ex}\log_2\!\left({\frac{T_n}{|V(\cG_n)|}}\right)\
\right)
\ee
Now let $t_n$ be the smallest integer with the property that
for all $v \in V(\cG_n)$, the subgraph of~$\cG_n$
induced by the neighborhood of $v$ has at most $t_n$ edges.
Then it follows from~\Eq{Bollobas} that
\be{JV-bound}
\alpha(\cG_n)
\ \ge \
\frac{(s_n+1)^n}{10 d_n}
\left( \Strut{4ex}{0ex}
\log_2 d_n \ - \
\shalf\hspace{0.25ex}\log_2\!\left({\frac{t_n}{3}}\right)\
\right)
\ee
Thus to obtain an asymptotic improvement over the Minkowski bound it would
suffice~to prove that $t_n=o(d_n^2)$. Before diving into the technical
details of the proof, let us explain intuitively why we expect to get 
an improvement by a factor that is linear in $n$.

Let us pick two points $x$ and $y$ uniformly at random in $\cS_n(2r_n)$.
The relevant question is: what is the probability that $d(x,y) < 2r_n$?
It is a rather standard fact in high-dimensional geometry that
(regardless of the value of $r_n$) this probability behaves as 
$e^{-cn}$ for large $n$. Therefore, we should expect that only 
an exponentially small fraction of pairs of points of~$\Z^n$ lying
within a sphere of radius $2r_n$ centered at some $z \in V(\cG_n)$ 
are adjacent in $\cG_n$. In other words, we expect that
$t_n \:\raisebox{2pt}{\footnotesize$\scriptstyle\lesssim$}\ d_n^2/e^{cn}$ 
which, in view of~(\ref{JV-bound}),
immediately leads to the desired $\Theta(n)$ improvement factor.
We derive a rigorous upper bound on $t_n$ next.

\looseness=-1
Consider the neighborhood of $\zero \in V(\cG_n)$, and let $\cH_n$
denote the subgraph of $\cG_n$ induced by this neighborhood.
As in Section\,2, we assume that the ratio $s_n/r_n$ in~\Eq{cubes} 
is sufficiently large so that 
$V(\cH_n) = (\Z^n {\setminus} \{\zero\}) \cap \cS_n(2r_n)$.
It is then obvious that $t_n = |E(\cH_n)|$, so 
\be{t_n-first}
2 t_n 
\ = \hspace{-2.5ex}
\sum_{\hspace*{2ex}x \inn V(\cH_n)} \hspace{-2.5ex} \deg(x)
\ee
where $\deg(x)$ denotes the degree of $x$ in $\cH_n$.
Write $S_1 = \cS_n(2r_n)$ and let $S_2$ be the sphere 
of radius $2r_n$ about $x \in V(\cH_n)$. Then $\deg(x)$
is just the number of points~of~$\Z^n$~in~$S_1 \cap S_2$.
Using the same argument as in~\Eq{d_n-upper}, we thus have
\be{degree-volume}
\deg(x) \ \le\ \vol(S'_1 \cap S'_2)
\ee
where $S'_1 = \cS_n(2r_n + \sqrt{n}/2)$ and $S'_2$ is the
sphere of radius $2r_n + \sqrt{n}/2$ about $x$. Clearly, 
the right-hand side of~\Eq{degree-volume} depends on $x$
only via its distance to the origin. Hence, define
\be{rho-delta}
\rho \ \ \deff\, \ \ 2r_n + \frac{\sqrt{n}}{2}
\hspace{4ex}\text{and}\hspace{4ex} 
\delta_x \ \ \deff\, \ \ \frac{d(x,\zero)}{2\rho}
\ee
with $\delta_x \in (0,\shalf)$ for all $x$.
It is not difficult to write down a precise expression for 
the~volume of $S'_1 \cap S'_2$ in terms of $\rho$ and~$\delta_x$.\hfill
Let $\theta = \cos^{-1} \delta_x$.\hfill Then $\vol(S'_1 \cap S'_2)$
is twice the\linebreak[3.99]\pagebreak[3.99] 

$\,$\vspace{1ex}

\begin{minipage}[t]{3.00in}
\centerline{\psfig{figure=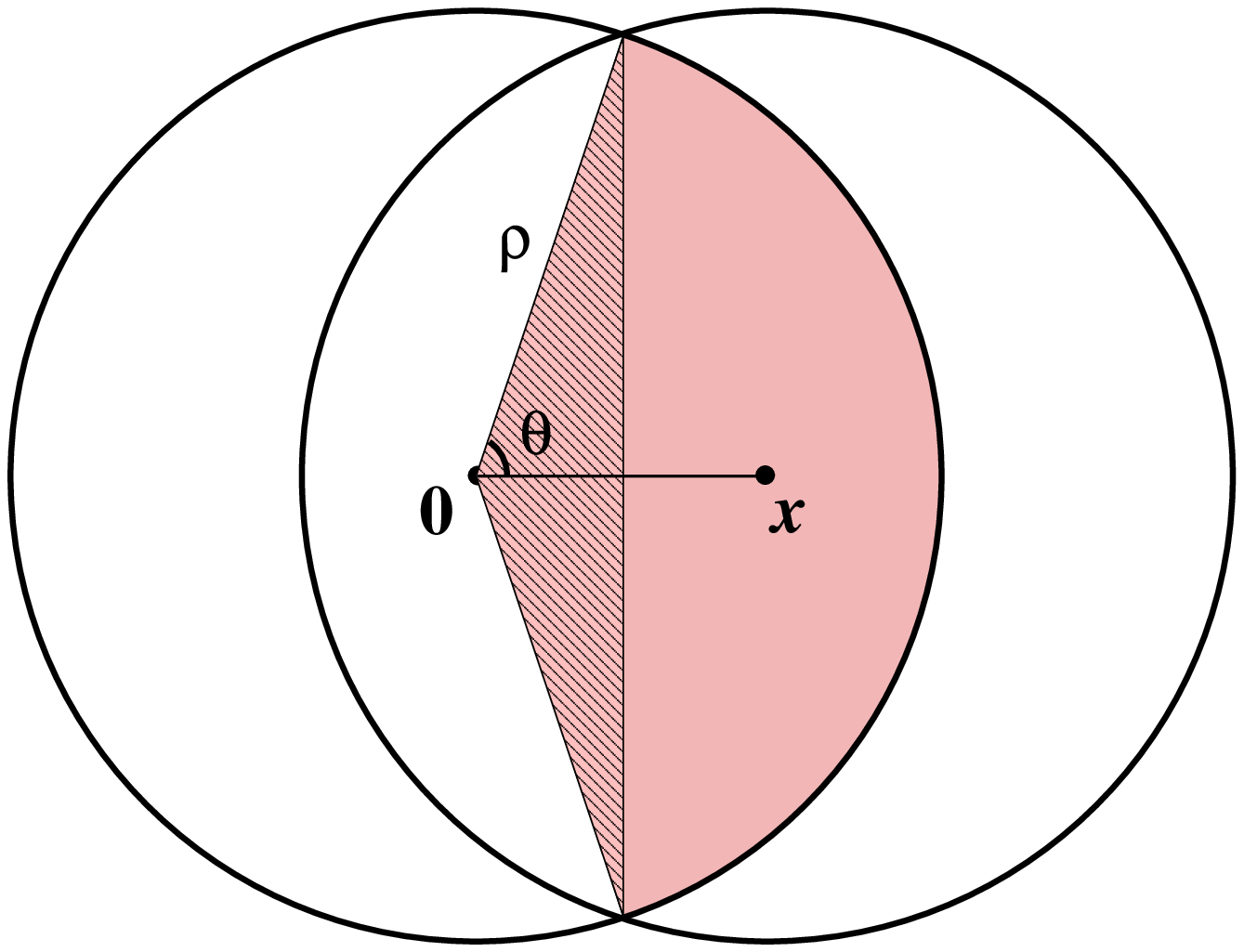,width=2.75in,silent=}}
\vspace{3.5ex}

\begin{center}
{\bf Figure\,1}  
\vspace{4.5ex}
\end{center}

\end{minipage}
\begin{minipage}[t]{3.00in}
\centerline{\psfig{figure=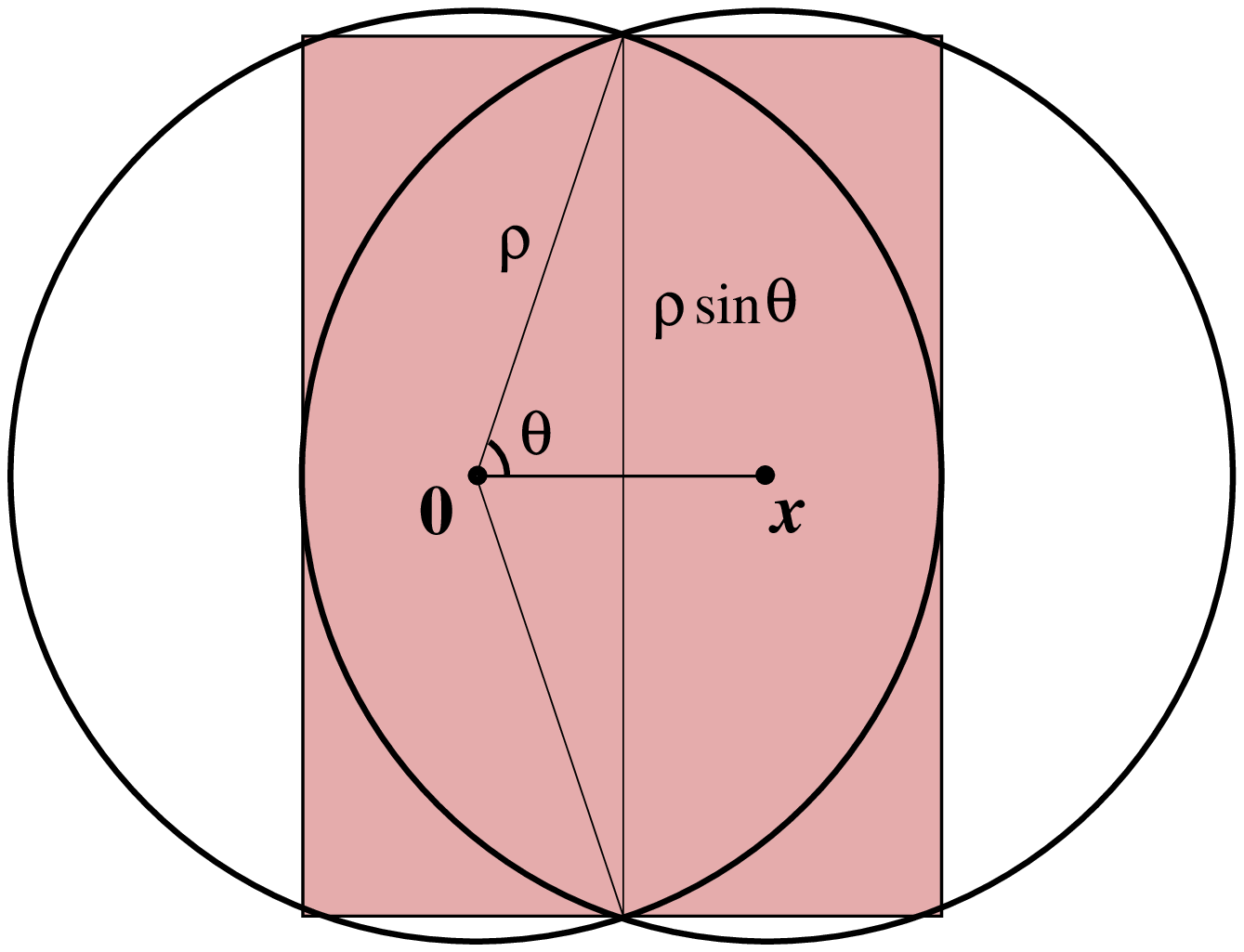,width=2.75in,silent=}}
\vspace{3.5ex}

\begin{center}
{\bf Figure\,2}  
\vspace{4.5ex}
\end{center}

\end{minipage}

\noindent
volume of a spherical sector of angle $2\theta$ 
(shaded in Fig.\,1) minus twice the volume of a~right cone 
of the same angle (cross-hatched in Fig.\,1). 
Thus
\be{exact-volume}
\vol(S'_1 \cap S'_2)
\ = \ 
\frac{2 \rho^n V_{n-1}}{n}
\left(
(n\,{-}\,1)\!\! \int_0^{\theta}\! (\sin\varphi)^{n-2} d\varphi
\,\ - \
\delta_x (\sin\theta)^{n-1} 
\right)
\ee
However, rather than estimating the integral 
$
\int_0^{\theta} (\sin\varphi)^{n-2} d\varphi
$
in \Eq{exact-volume}, we will use the following simple bound (without
compromising much in the asymptotic quality of the obtained result).
It is easy to see (cf.~Fig.\,2)
that $S'_1 \cap S'_2$ is contained in a cylinder~of height
$2\rho - d(x,\zero)$
whose base is an $(n{-}1)$-dimensional sphere of radius $\rho\sin\theta$.
Hence
\be{cylinder-volume}
\vol(S'_1 \cap S'_2)
\ \le \ 
2\rho\,(1\,{-}\,\delta_x) V_{n-1} \rho^{n-1} (\sin\theta)^{n-1}
\ \le \
(1\,{-}\,\delta^2_x)^{n/2}\, n V_n \rho^n
\ee
where the second inequality follows from the fact that 
$2V_{n-1} \le nV_n$ for all $n$.
Now let $U_k$ denote the set of all $x \in V(\cH_n)$ such that
$d^2(x,\zero) = k$. Clearly $d^2(x,\zero)$ is an integer in the range
$1 \le d^2(x,\zero) \le 4r_n^2 - 1$. We thus break the sum 
in~\Eq{t_n-first} into two parts
\be{t_n-second}
2 t_n 
\ = \hspace{-2.5ex}
\sum_{\hspace*{2ex}x \inn V(\cH_n)} \hspace{-2.5ex} \deg(x)
\:\ = \,
\sum_{k=1}^{\sfourth r^2_n-1} \hspace{-1.15ex}
\sum_{x \inn U_k} \hspace{-.25ex} \deg(x)
\:\ + \,
\sum_{k={\sfourth r^2_n}}^{4r^2_n-1} 
\sum_{x \inn U_k} \hspace{-.25ex} \deg(x)
\ee
and bound each part separately. 
As it turns out,
crude upper bounds~on~$|U_k|$~suffice in each case. For the
first double-sum in~\Eq{t_n-second}, we use the fact that
\mbox{$\deg(x) \le d_n \le V_n\rho^n$}
for all $x \in V(\cH_n)$ by~\Eq{d_n-bounds}. Therefore, 
applying once again the method of~\Eq{d_n-upper}, we have
\be{first-sum}
\hspace*{-1.0ex}\sum_{k=1}^{\frac{r^2_n}{4}-1} \hspace{-.50ex}
\sum_{x \inn U_k} \hspace{-.5ex} \deg(x)
\ \le \
V_n \rho^n \hspace{-.4ex} \sum_{k=1}^{\frac{r^2_n}{4}-1} \hspace{-.5ex} |U_k|
\ \le \
V_n \rho^n\hspace{0.40ex}
\vol\!\left(\hspace{-.2ex}\Strut{4.00ex}{0ex}
\cS_n\!\left(\hspace{-.2ex}
\frac{r_n \,{+} \sqrt{n}}{2}\right)\hspace{-.6ex}\right)
\ \le \
\left(\frac{1}{2}\right)^{\hspace{-.50ex}n} \hspace{-.7ex} 
V_n^2 \rho^{2n} 
\ee
where the last inequality assumes $r_n \ge \sqrt{n}/2$,
so that $\rho \ge r_n + \sqrt{n}$. In fact, henceforth,~let 
us take $r_n = 2n^2$ as in Section\,2. Then $k \ge n^4$ in the second 
sum of~\Eq{t_n-second},\pagebreak[3.99]
and we can bound $|U_k|$ as follows:
$
|U_k|
\le
\vol\bigl(\cS_n(\sqrt{k}+\sqrt{n}/2)\bigr)
\le
V_n k^{n/2} \bigl(1 + \sqrt{n}/(2n^2)\bigr)^n
\le 
2 V_n k^{n/2} 
$.
Combining this with the bounds \Eq{degree-volume} and
\Eq{cylinder-volume} on $\deg(x)$, we have
$$
\sum_{k={\sfourth r^2_n}}^{4r^2_n-1} 
\sum_{x \inn U_k} \hspace{-.35ex} \deg(x)
\ \le \
n V_n \rho^n \hspace{-.5ex} 
\sum_{k={\sfourth r^2_n}}^{4r^2_n-1} \hspace{-.75ex} 
|U_k| \hspace{-.25ex} \left(1 - \delta_k^2\right)^{\hspace{-.25ex}n/2}
\hspace{-.5ex}\le \
2^{n+1} n V^2_n \rho^{2n} \hspace{-.5ex} 
\sum_{k={\sfourth r^2_n}}^{4r^2_n-1} \hspace{-.5ex} 
\left(\delta_k^2\, (1 {-} \delta_k^2)\right)^{\hspace{-.25ex}n/2}
$$
where $\delta_k \ \deff\  \sqrt{k}/(2\rho)$. Now, the function 
$f(\delta) = \delta^2(1 - \delta^2)$ attains its maximum 
in the interval $[0,\shalf]$ at $\delta = \shalf$. 
Putting this together with \Eq{t_n-second} and \Eq{first-sum},
we finally obtain the desired upper bound on $t_n$, namely
\be{t_n-bound}
t_n
\ \le \
\left(\frac{1}{2}\right)^{\hspace{-.35ex}n+1} \hspace{-1.0ex} 
V_n^2 \rho^{2n} 
\hspace{1.15ex} + \hspace{1.35ex} 
2^{n} n V^2_n \rho^{2n} \hspace{-.5ex} 
\sum_{k={\sfourth r^2_n}}^{4r^2_n-1} \hspace{-.5ex} 
\left(\frac{3}{16}\right)^{\hspace{-.35ex}n/2}
\hspace{-0.25ex} \le \
\left(\hspace{-0.35ex}\frac{\sqrt{3}}{2}\right)^{\hspace{-.35ex}n}
\hspace{-.25ex}n V_n^2\,\rho^{2n+2}
\ee
Recall that $d_n\le V_n\rho^n$ by~\Eq{d_n-bounds}.
Substituting this bound together with~\Eq{t_n-bound} in the
bound~\Eq{JV-bound} on the independence number of $\cG_n$
produces
\begin{eqnarray*}
\alpha(\cG_n)
&\hspace{-.5ex}\ge\hspace{-.25ex}& 
\frac{(s_n\,{+}\,1)^n}{10 V_{\!n}\hspace{.10ex}\rho^n}
\left(\Strut{5ex}{0ex}
\log_2\hspace{-.5ex}\left(\Strut{2ex}{0ex}V_{\!n}\hspace{.10ex}\rho^n\right)
\ - \
\frac{1}{2}\log_2\!\left(\hspace{-.75ex}
\left(\hspace{-.50ex}\frac{\sqrt{3}}{2}\right)^{\hspace{-.40ex}n\phantom{1}}
\hspace{-1.50ex}V_n^2\rho^{2n}\right)
\ - \
\frac{1}{2}\log_2\hspace{-.25ex}\left(n\rho^2\right)
\right)\\[1ex]
&\hspace{-.5ex}=\hspace{-.25ex}& 
\frac{(s_n\,{+}\,1)^n}{10\, V_{\!n} \left(2r_n + \sqrt{n}/2\right)^n}
\left(
\frac{n\,\log_2\bigl(2/\!\sqrt{3}\,\bigr)}{2}
\:\ - \:\
\frac{\log_2\bigl(n(2r_n+\sqrt{n}/2)^2\,\bigr)}{2}
\right)
\end{eqnarray*}
where we have used the definition of $\rho$ in~\Eq{rho-delta}.
Finally, using \Eq{first-bound} with 
$|\cI| = \alpha(\cG_n)$ while taking $r_n = 2n^2$ and $s_n = 2n^4$ as before,
we obtain the following bound 
\begin{eqnarray}
\Delta(\cP_n)
&\hspace{-.85ex}\ge\hspace{-.50ex}& 
\frac{n}{20} \hspace{-.25ex}
\left(\hspace{-1pt} \frac{s_n}{s_n\,{+}\,2r_n}
\hspace{-1pt}\right)^{\hspace{-.50ex}n}\hspace{-.55ex} 
\left(\hspace{-1pt}\frac{r_n}{2r_n\,{+}\sqrt{n}/2}
\hspace{-1pt}\right)^{\hspace{-.5ex}n}\hspace{-.20ex}
\left(\hspace{-.25ex}
\log_2\bigl(2/\!\sqrt{3}\,\bigr)
\, - \,
\frac{\log_2\bigl(n(2r_n\,{+}\sqrt{n}/2)^2\bigr)}{n}
\right)\hspace*{-5pt}\nonumber\\[1ex]
&\hspace{-.85ex}=\hspace{-.50ex}& 
\frac{\log_2\bigl(2/\!\sqrt{3}\,\bigr)}{20}\
n 2^{-n} \left(\Strut{2.5ex}{0ex} 1 + o(1)\right)
\label{final-bound}
\end{eqnarray}
Since $\log_2\bigl(2/\!\sqrt{3}\,\bigr)/20 = 0.0103..$,
this establishes~\Eq{thm1}. 
Now, given any graph~\mbox{$G = (V,E)$},
there is a deterministic algorithm~\cite{HL} that finds 
an independent set $\cI$ in $G$ whose size is lower bounded
by~\Eq{Bollobas} in time $O(d_{\rm av}|E| + |V|)$, where 
$d_{\rm av}$ is the average degree of $G$. In the case
of the graph $\cG_n$, this reduces to
$
O\!\left(V_n^2 (s_n+1)^n (2r_n \,{+} \sqrt{n}/2)^{2n} \right)
$.
With~\mbox{$r_n = 2n^2$}~and $s_n = 2n^4$, the expression
$
V_n^2 (s_n+1)^n (2r_n + \sqrt{n}/2)^{2n} 
$
behaves as $(64\pi e n^7)^n$ for large~$n$.
This completes the proof of \Tref{Thm1}, 
with the value of $\gamma$ given by $7 + \epsilon$.
However,~note that our choice of the values $r_n = 2n^2$ and $s_n = 2n^4$
was motivated primarily by~the notational convenience of having
$r_n,s_n \in 2\Z$. In fact 
$r_n = n^{1.5 + \epsilon}$ and $s_n = n^{2.5 + \varepsilon}$ 
would suffice, \looseness=-1
as can be readily seen from~\Eq{Minkowski-bound} and~\Eq{final-bound}. 
Thus the value of $\gamma$ can be taken as $4.5+\epsilon$.


\end{document}